\documentclass[11pt]{article}
\usepackage{amsmath}
\usepackage{graphics, amssymb}%, showlabels}
\newtheorem{Def}{Definition}[section]
\newtheorem{Thm}[Def]{Theorem}
\newtheorem{Prp}[Def]{Proposition}
\newtheorem{Lemma}[Def]{Lemma}

\newtheorem{Corollary}[Def]{Corollary}

\newcommand{\Proof}{{\em{Proof. }}}
\newcommand{\QED}{\ \hfill $\FBox$ \\[1em]}
\newcommand{\spc}{\;\;\;\;\;\;\;\;\;\;}
\newcommand{\R}{\mathbb{R}}
\newcommand{\C}{\mathbb{C}}

\newcommand{\FBox}{\rule{2mm}{2.25mm}}
\newcommand{\ov}{\overline}

\newcommand{\D}{\mathcal{D}}

\newcommand{\spec}{\mbox{\rm spec}}

\newcommand{\Tr}{\mbox{\rm{Tr}\/}}
\newcommand{\beq}{\begin{equation}}
\newcommand{\eeq}{\end{equation}}

\numberwithin{equation}{section}

\title{A Level Set Analysis of the Witten Spinor \\
with Applications to Curvature Estimates}
\author{Felix Finster\thanks{Supported by the Deutsche Forschungsgemeinschaft
within the Priority Program ``Globale Differentialgeometrie''.}}
\date{February 2007}

\begin{document}
\maketitle
\begin{abstract}
We analyze the level sets of the norm of the Witten spinor in an
asymptotically flat Riemannian spin manifold of positive scalar curvature.
Level sets of small area are constructed. We prove curvature estimates
which quantify that, if the total mass becomes small, the
manifold becomes flat with the exception of a set of small surface area.
These estimates involve either a volume bound or a spectral bound
for the Dirac operator on a conformal compactification, but they are independent
of the isoperimetric constant.
\end{abstract}

\section{Introduction and Statement of Results} \label{sec1}
Asymptotically flat manifolds of positive scalar curvature describe
isolated time-symmetric gravitating systems in general relativity.
The main point of mathematical interest is to explore the connections
between the total mass and the local geometry of the manifold.
A common general method for analyzing related questions is to
consider a flow of hypersurfaces, under which a certain quasi-local
mass functional is monotone. The most prominent example is the
inverse mean curvature flow, under which the Hawking mass is
monotone~\cite{HI}. In order to conveniently parametrize the
hypersurfaces, one often represents the hypersurfaces as level sets
of a real-valued function~$\phi$, which is then a solution of a suitable
partial differential equation on the manifold.
In this paper, we follow this approach, taking for the function~$\phi$
the norm~$|\psi|$ of a Witten spinor in an asymptotically flat spin
manifold~$M$.
This is particularly simple because the well-known
existence of the Witten spinors ensures global existence of the corresponding
flow. Nevertheless, this flow has nice and useful properties,
above all that, in analogy to a monotonicity property, integrating~$|D\phi|^2$ over the set~$\{x \in M \:|\: \phi(x) < \tau\}$ gives a {\em{convex}}
function in~$\tau$.

As an application we construct level sets $\{x \:|\: \phi(x)=t\}$
of small area. Combining this result with the curvature
estimates~\cite{BF, FK}, we prove estimates of the following type:
There is an exceptional set~$\Omega \subset M$ of small surface
area~$|\partial \Omega|$ such that on~$M \setminus \Omega$ the
Riemann tensor is small in an $L^2$-sense. Here by ``small'' we mean that
the upper bounds involve positive powers of the total mass~$m$, and thus
tend to zero as~$m \searrow 0$.
These curvature estimates are stated in two versions, either involving
a volume bound, or, using methods and results from~\cite{FKr},
involving a spectral bound for the Dirac operator on a conformal
compactification of~$M$.

We now introduce the mathematical framework and state our results.
Let~$(M^n, g)$ be a complete Riemannian spin manifold of dimension $n \geq 3$.
We assume that the {\em{scalar curvature}} of~$M$ is
{\em{non-negative}} and integrable. Furthermore,
we assume that~$M$ is {\em{asymptotically flat}}, for simplicity
with one asymptotic end. Thus there is a compact set~$K \subset M$
and a diffeomorphism $\Phi \,:\, M \setminus K \rightarrow
\R^n \setminus B_{\rho}(0)$, $\rho>0$, such that
\[ (\Phi_* g)_{ij} \;=\; \delta_{ij} \:+\: {\cal{O}}(r^{2-n}) \:,\quad
\partial_k (\Phi_* g)_{ij} \;=\; {\cal{O}}(r^{1-n}) \:,\quad
\partial_{kl} (\Phi_* g)_{ij} \;=\; {\cal{O}}(r^{-n}) \:. \]
Under these assumptions, the {\em{total mass}} of the manifold is
defined by~\cite{ADM}
\begin{equation}
m \;=\; \frac{1}{c(n)} \lim_{\rho \rightarrow \infty}
\int_{S_\rho} (\partial_j (\Phi_* g)_{ij} - \partial_i (\Phi_* g)_{jj})
\:d\Omega^i \:, \label{massdef}
\end{equation}
where $c(n)>0$ is a normalization constant
and $d\Omega^i$ denotes the product of the volume form
on~ $S_\rho \subset \R^n$ by the $i$-th
component of the normal vector on $S_\rho$.

Spinors are very useful for the analysis of asymptotically flat
spin manifolds. The basic identity is the
{\em{Lichnerowicz-Weitzenb{\"o}ck formula}}
\beq \label{LW}
\D^2 \;=\; -\nabla^2 + \frac{s}{4} \;,
\eeq
where~$\D$ is the Dirac operator, $\nabla$ is the spin connection,
and~$s$ denotes scalar curvature.
Witten~\cite{W} considered solutions of the Dirac equation with constant
boundary values~$\psi_0$ in the asymptotic end,
\beq \label{bvp}
\D \psi \;=\; 0 \:,\qquad
\lim_{|x| \rightarrow \infty} \psi(x) \;=\; \psi_0\:,
\eeq
where~$\psi$ is a smooth section of the spinor bundle~$SM$.
In~\cite{PT, B} it is proven that for any~$\psi_0$,
this boundary value problem has a unique solution. We refer to~$\psi$ as the
{\em{Witten spinor}} with boundary values~$\psi_0$.
For a Witten spinor, the Lichnerowicz-Weitzenb\"ock formula implies that
\beq \label{div}
\nabla_i \,\langle \psi, \nabla^i \psi \rangle \;=\;
|\nabla \psi|^2 + \frac{s}{4}\: |\psi|^2\:.
\eeq
Integrating over~$M$, applying Gauss' theorem and relating the boundary
values at infinity to the total mass
(where we choose~$c(n)$ in~(\ref{massdef}) appropriately),
one obtains the identity~\cite{W, PT, B}
\beq \label{ibp}
\int_M \left( |\nabla \psi|^2 + \frac{s}{4}\, |\psi|^2 \right) d\mu_M \;=\; m\:.
\eeq
This identity immediately implies the
{\em{positive mass theorem}} for spin manifolds
(for the positive mass theorem on non-spin manifolds see~\cite{SY, L}).

In this paper, we consider the {\em{level sets}} of the norm of a
Witten spinor~$\psi$. We introduce the function $\phi=|\psi|$
and set
\[ \tau_0 \;=\; \inf_M \phi \:,\qquad \tau_1 \;=\; \sup_M \phi\:. \]
For any~$\tau \in [\tau_0, \tau_1]$ we define the set
\[ \Omega(\tau) \;=\; \left\{ x \in M \:|\:
\phi(x) < \tau \right\} \:. \]
Clearly, the sets~$\Omega(\tau)$ are open and form an increasing
family in the sense that~$\tau' \leq \tau$
implies~$\Omega(\tau) \subset \Omega(\tau')$.
Moreover, the boundary of~$\Omega(\tau)$
is the level set~$\{ x \:|\: \phi(x)= \tau \}$.
We also introduce the two functions
\begin{eqnarray}
m(\tau) &=& \int_{\Omega(\tau)}
\left( |\nabla \psi|^2 + \frac{s}{4}\, |\psi|^2 \right) d\mu_M \label{mudef} \\
F(\tau) &=& \int_{\Omega(\tau)} |D \phi|^2\, d\mu_M\:. \label{Fdef}
\end{eqnarray}
Since the integrands are non-negative, it is obvious that these
functions are monotone increasing.
Furthermore, $m(\tau_0)=0=F(\tau_0)$.
Comparing~(\ref{mudef}) with~(\ref{ibp}), one sees that
\[ m(\tau_1) \;=\; m \:; \]
this is why we refer to~$m(\tau)$ as the {\em{mass function}}.
Furthermore, the Kato inequality
$|D |\psi|| \leq |\nabla \psi|$ implies that
\beq \label{Fmineq}
F(\tau) \;\leq\; m(\tau)\:.
\eeq
Our main result relates~$F'$ to~$m$ and makes a convexity statement.
\begin{Thm} \label{thm1}
The function $F(\tau) : [\tau_0, \tau_1] \rightarrow \R$
is convex and differentiable almost everywhere.
It satisfies for almost all~$\tau$ the identity
\beq \label{deq}
\frac{dF(\tau)}{d\tau} \;=\; \frac{1}{\tau}\: m(\tau) \:.
\eeq
\end{Thm}
As an immediate application, this theorem implies the inequalities
\[ \frac{m(\tau')}{\tau'} \;\leq\; \frac{F(\tau) - F(\tau')}{\tau-\tau'}
\;\leq\; \frac{m(\tau)}{\tau} \spc
{\mbox{for almost
all~$\tau, \tau' \in [\tau_0, \tau_1]$ with~$\tau' < \tau$}}\:, \]
giving useful information on the behavior of the Witten spinor.
For example, setting~$\tau'=\tau_0$, the right inequality
yields the following upper bound for~$\tau_0$,
\[ \tau_0 \;\leq\; \tau \left( 1 - \frac{F(\tau)}{m(\tau)} \right) . \]
Inequalities of this type seem surprising. However, one should keep in
mind that in all interesting applications, the function~$F$ is difficult
to compute, and therefore these inequalities are of limited practical value.

Here we focus on applications of the above theorem to
{\em{curvature estimates}} in asymptotically flat manifolds
in the spirit of~\cite{BF, FK}. The main point is that we
now obtain estimates which do not depend on the isoperimetric
constant, and where the exceptional set has small surface area
(instead of small volume as in~\cite{BF, FK}).
Here we state the results in the physically interesting case
of dimension~$n=3$, but we also prove similar results in
general dimension (see Theorems~\ref{thmg1} and~\ref{thmg2} below).
We begin with a curvature estimate assuming a volume bound for some
set~$\Omega(t_1) \setminus \Omega(t_0)$.
\begin{Thm} \label{thm31}
Let~$(M^3, g)$ be a complete, asymptotically flat manifold whose
scalar curvature is non-negative and integrable.
For any Witten spinor~$\psi$ (\ref{bvp}) and any interval~$[t_0, t_1]
\subset (0,1]$, there is $t \in [t_0, t_1]$
with the following properties. The level set~$|\psi|=t$ is
a submanifold of~$M$, whose $2$-volume $A(t)$ is bounded by
\beq \label{Aineq}
A(t) \;\leq\; \sqrt{F(t_1)-F(t_0)}\; \frac{\sqrt{V(t_1) - V(t_0)}}{t_1-t_0} \:,
\eeq
where~$V(t) := \mu_M(\Omega(t))$.
On the set~$M \setminus \Omega(t)$, the Riemann tensor satisfies the inequality
\[ \int_{M \setminus \Omega(t)} |R|^2 \;\leq\;
\frac{m \:c_1}{t^2}\: \sup_M |R| \:+\:
\frac{\sqrt{m}\:c_2}{t^2}\: \|\nabla R\|_{L^2(M)} \]
with constants~$c_1, c_2$ which are independent of the geometry of~$M$.
\end{Thm}
Note that by~(\ref{Fmineq}), we can always bound~$F$ by the total mass.
Furthermore, it is in most applications sufficient to drop the
term~$V(t_0)$ in~(\ref{Aineq}) and to choose~$t_0=t_1/2$.
This gives the following corollary.
\begin{Corollary} \label{coroll}
Let~$(M^3, g)$ be a complete, asymptotically flat manifold 
whose scalar curvature is non-negative and integrable.
For any Witten spinor~$\psi$ (\ref{bvp}) and any~$t_1 \in (0,1]$, there is $t
\in [\frac{t_1}{2}, t_1]$ with the following properties. The level set~$|\psi|=t$ is
a submanifold of~$M$ whose $2$-volume $A(t)$ is bounded by
\[ A(t) \;\leq\; \frac{2\,\sqrt{m\, V(t_1)}}{t_1} \:. \]
On the set~$M \setminus \Omega(t)$, the Riemann tensor satisfies the
inequality
\[ \int_{M \setminus \Omega(t)} |R|^2 \;\leq\;
\frac{m \:c_1}{t_1^2}\: \sup_M |R| \:+\:
\frac{\sqrt{m}\:c_2}{t_1^2}\: \|\nabla R\|_{L^2(M)} \;. \]
\end{Corollary}

The remaining question is how to control the volume~$V(t_1) - V(t_0)$.
We here propose the method to work with a {\em{spectral bound}}
for the Dirac operator {\em{on a conformal compactification}} of~$M$.
As in~\cite{FKr} we assume for simplicity that~$M$ is
{\em{asymptotically Schwarzschild}}, although the method should apply
to more general asymptotically flat manifolds as well.
Thus under the diffeomorphism~$\Phi \,:\, M \setminus K \rightarrow \R^n \setminus B_{\rho}(0)$ the metric becomes the
Schwarz\-schild metric,
\beq \label{assch}
(\Phi_* g)_{ij} \;=\; \left(1+\frac{m}{2r} \right)^{4}\:\delta_{ij}
\:.
\eeq
We point compactify~$M$ by a conformal transformation
\[ \tilde{g} \;=\; \lambda^2\: g \:, \]
in such a way that the geometry of~$K$ remains unchanged,
the scalar curvature stays non-negative,
and the compactification of the asymptotic end is isometric
to a cap~$C_R \subset S^n_\sigma$ of a sphere of radius~$\sigma$
(for details see~\cite{FKr} and Section~\ref{sec4}).
The Dirac operator on the compactification~$(\bar{M}, \tilde{g})$
is denoted by~$\tilde{\D}$.
\begin{Thm} \label{thm32}
Let~$(M^3, g)$ be a complete, manifold of non-negative
scalar curvature such that~$M \setminus K$ is isometric to
the Schwarzschild geometry~(\ref{assch}).
Then there is~$t \in (\frac{1}{4},\frac{1}{2})$
such that the level set~$|\psi|=t$ is
a submanifold of~$M$. Its $2$-volume $A(t)$ is bounded by
\[ A(t) \;\leq\; c_0\, \sqrt{m}
\;\frac{(\rho+m)^{\frac{3}{2}}}{\sigma\, \inf \spec |\tilde{\D}|}\: . \]
On the set~$M \setminus \Omega(t)$, the Riemann tensor satisfies the inequality
\[ \int_{M \setminus \Omega(t)} |R|^2 \;\leq\;
m \:c_1\: \sup_M |R| \:+\:
\sqrt{m}\:c_2\: \|\nabla R\|_{L^2(M)} \]
with constants~$c_0$, $c_1$ and~$c_2$, which are independent of the
geometry of~$M$.
\end{Thm}

\section{Level Set Analysis} \label{sec2}
Let~$\psi$ be a Witten spinor~(\ref{bvp}). By linearity,
it is no loss of generality to always normalize~$\psi_0$ by one,
\beq \label{norm}
|\psi_0| \;=\; 1\:.
\eeq
We introduce the level sets of~$\phi$ by
\[ L(\tau) \;=\; \left\{ x \in M \:|\: \phi(x) = \tau \right\} . \]
We call~$\tau \in [\tau_0, \tau_1]$ a {\em{regular value}}
if~$D \phi$ has no zeros on~$L(\tau)$, otherwise it is called a
{\em{singular value}}. According to Sard's lemma, the singular
values form a set of Lebesgue measure zero in~$[\tau_0, \tau_1]$.
If~$\tau$ is a regular value, the implicit function theorem yields
that~$L(\tau)$ is a smooth submanifold of~$M$ of codimension one.
In this case, we denote the induced measure on~$L(\tau)$
by~$d\mu_{L(\tau)}$. 

We first motivate our method, neglecting the subtle issue of the
singular values.
A promising idea for getting information on the level sets is to
integrate a smooth function~$h$ on~$M$ (which may be a curvature
expression or an expression involving the Witten spinor) over the
level sets,
\[ \int_{L(\tau)} h(x)\: d\mu_{L(\tau)}(x) \:, \]
and to analyze a ``flow equation'' for this expression.
To derive the flow equation, we first apply Gauss' theorem to obtain
\[ \int_{L(\tau)} h\, d\mu_{L(\tau)} \;=\; \int_{\Omega(\tau)}
\nabla_i \left( h\, \nu^i \right) d\mu_M \:, \]
where~$\nu^i = (D^i \phi)/|D\phi|$ denotes the outer normal.
Now the co-area formula yields
\[ \int_{L(\tau)} h\, d\mu_{L(\tau)} \;=\;
\int_{\tau_0}^\tau d\sigma \int_{L(\sigma)} \frac{1}{|D\phi|}\:
\nabla_i \left( \frac{h\: D^i \phi}{|D \phi|} \right) d\mu_{L(\sigma)}\:, \]
and differentiating with respect to~$\tau$ gives for any regular
value~$\tau$ the differential equation
\beq \label{floweq}
\frac{d}{d\tau} \int_{L(\tau)} h\, d\mu_{L(\tau)} \;=\;
\int_{L(\tau)} \frac{1}{|D\phi|}\:
\nabla_i \left( \frac{h\: D^i \phi}{|D \phi|} \right) d\mu_{L(\tau)}\:.
\eeq
The basic problem is that the right hand side will in general
involve new geometric quantities, which are difficult to control.
For example, setting~$h \equiv 1$, we obtain the area
functional~$A(\tau):=\mu_{L(\tau)}(L(\tau))$.
Its flow equation is
\[ \frac{d}{d\tau} A(\tau) \;=\;
\int_{L(\tau)} \frac{1}{|D\phi|}\:
\nabla_i \left( \frac{D^i \phi}{|D \phi|} \right) d\mu_{L(\tau)}\:. \]
Here in the integrand the well-known mean curvature operator appears.
But the mean curvature of the level sets is not known, and it seems
difficult to get information on mean curvature. Therefore,
in order to make use of~(\ref{floweq}), we must look for a special
function~$h$ for which the right side of~(\ref{floweq}) has nice properties.
Choosing~$h=|D \phi|$, we get the simple equation
\[ \frac{d}{d\tau} \int_{L(\tau)} |D\phi|\, d\mu_{L(\tau)} \;=\;
\int_{L(\tau)} \frac{\Delta \phi}{|D\phi|}\:d\mu_{L(\tau)}\:. \]
The integral on the left equals the function~$f$ which we shall define
below. It is preferable to introduce it using volume integrals
over~$\Omega(\tau)$ instead of surface integrals, because such volume integrals
make sense even if~$\tau$ is a singular value. This motivates the following constructions.

A direct calculation using the
Licherowicz-Weitzenb\"ock formula~(\ref{LW}) gives
\begin{eqnarray*}
\Delta \phi^2 &=& \Delta \langle \psi, \psi \rangle \;=\;
2 \,{\mbox{Re}} \langle \nabla^2 \psi, \psi \rangle \:+\: 2 \,|\nabla \psi|^2
\;=\; \frac{s}{2}\: |\psi|^2 + 2 \,|\nabla \psi|^2 \nonumber \\
\nabla \phi &=& \frac{\nabla \phi^2}{2 \phi} \nonumber \\
\Delta \phi &=& \frac{\Delta \phi^2}{2 \phi} - \frac{(\nabla \phi^2)^2}{4 \phi^2} \;=\; \frac{s}{4}\, \phi \:+\: \frac{|\nabla \psi|^2}{\phi}
\:-\: \frac{|{\mbox{Re}} \langle \nabla \psi, \psi \rangle|^2}{\phi^3} \:.
\end{eqnarray*}
Applying the the Schwarz inequality~$|{\mbox{Re}} \langle \nabla \psi, \psi \rangle| \leq |\nabla \psi| \, \phi$,
we obtain the inequality
\beq \Delta \phi \;\geq\; \frac{s}{4}\: \phi \:. \label{sub} \eeq
It is an important observation that~$\phi$ is subharmonic.
In particular, we can apply the maximum principle to conclude
that~$|\phi| \leq |\psi_0|$. Comparing with~(\ref{norm}), we find that
\beq \label{22a}
\tau_1 \;=\; 1 \:.
\eeq

We introduce the function~$f(\tau)$ by
\beq \label{fdef}
f(\tau) \;=\; \int_{\Omega(\tau)} \Delta \phi\: d\mu_M\:.
\eeq

\begin{Lemma} \label{lemma1}
The function $f$ is monotone increasing and left-sided continuous,
i.e.\ for all~$\tau \in (\tau_0, 1]$,
\[ \lim_{\tau' \nearrow \tau} f(\tau') \;=\; f(\tau)\:. \]
For almost all~$\tau \in [\tau_0,1]$,
\beq \label{fob}
f(\tau) \;=\; \int_{L(\tau)} |D \phi|\: d\mu_{L(\tau)}\:.
\eeq
\end{Lemma}
\Proof We write~$f$ in the form
\[ f(\tau) \;=\; \int_M g\, d\mu_M \quad {\mbox{with}} \quad
g(x) \;:=\; \Delta \phi(x)\; \chi_{\Omega(\tau)}(x)\:, \]
where~$\chi$ denotes the characteristic function.
According to~(\ref{sub}), the integrand is non-negative. The monotonicity
of~$f$ is obvious because the family~$\Omega(\tau)$ is increasing.
To prove left-sided continuity, we note that for all~$\tau' < \tau$
\[ f(\tau) - f(\tau') \;=\; \int_M \Delta \phi(x)\; \chi_{\Omega(\tau)
\setminus \Omega(\tau')}(x)\: d\mu_M\:. \]
As~$\tau' \nearrow \tau$, the characteristic function tends to zero
pointwise. Hence in this limit, $f(\tau) - f(\tau')$ tends to zero
due to Lebesgue's monotone convergence theorem. 

To prove~(\ref{fob}), we let~$\tau$ be a regular value of~$\phi$.
Then the outer normal on~$L(\tau)$ is given by
\[ \nu \;=\; \frac{D \phi}{|D \phi|}\:. \]
Thus applying Gauss' theorem in~(\ref{fdef}), we obtain
\[ f(\tau) \;=\; \int_{\Omega(\tau)} \nabla_i \left( D^i \phi \right)
\;=\; \int_{L(\tau)} (D^i \phi)\: \nu_i\, d\mu_{L(\tau)} \;=\;
\int_{L(\tau)} |D \phi|\, d\mu_{L(\tau)} \:. \]

\vspace*{-0.6cm}
\QED

\vspace*{.3em} \noindent
{\em{Proof of Theorem~\ref{thm1}. }}
The co-area formula yields
\beq \label{coarea}
F(\tau) \;=\; \int_{\tau_0}^\tau \left[ \int_{L(\sigma)}
|D \phi|\: d\mu_{L(\sigma)} \right] d\sigma\:,
\eeq
where the square bracket is defined almost everywhere according
to Sard's lemma. Applying Lemma~\ref{lemma1}, the square
bracket coincides with~$f$ and is thus monotone increasing and
left-sided continuous. This implies that~$F$ is~$C^0([\tau_0, 1])$,
and is differentiable almost everywhere with~$F'(\tau)=f(\tau)$.
Since~$f$ is monotone increasing, we conclude that~$F$ is convex.

It remains to show that for almost all~$\tau \in [\tau_0, 1]$,
\[ f(\tau) \;=\; \frac{1}{\tau}\: m(\tau)\:. \]
According to Sard's lemma, we may assume that~$\tau$ is a regular
value. Then, applying~(\ref{div}) and Gauss' theorem in~(\ref{mudef}),
we obtain
\[ m(\tau) \;=\; \int_{\Omega(\tau)} {\mbox{Re}}\,
\nabla_i \,\langle \psi, \nabla^i \psi \rangle\: d\mu_M 
\;=\; \int_{L(\tau)} {\mbox{Re}}\,
\langle \psi, \nabla^\nu \psi \rangle\: d\mu_{L(\tau)} \:. \]
Using furthermore that $D \phi = {\mbox{Re}} \langle \psi, \nabla \psi \rangle/\phi$, and that~$D \phi$ points in normal direction, we obtain
\[ m(\tau) \;=\; \int_{L(\tau)} \phi\, D_\nu \phi\: d\mu_{L(\tau)}
\;=\; \tau \int_{L(\tau)} |D \phi|\: d\mu_{L(\tau)}
\;=\; \tau\, f(\tau)\:. \]

\vspace*{-0.6cm}
\QED

Next we want to derive an inequality involving the volume of the sets~$\Omega(\tau)$ and the area of the level sets. We define the volume~$V(\tau)$ by
\[ V(\tau) \;=\; \int_{\Omega(\tau)} d\mu_M \:,\qquad \tau \in [\tau_0, 1]\:. \]
The area of the level sets is only defined if~$\tau$ is a regular value,
and we extend the area function by zero to the singular values,
\[ A(\tau) \;=\; \left\{ \begin{array}{cl}
\mu_{L(\tau)}(L(\tau)) & {\mbox{if~$\tau$ is a regular value}} \\[.8em]
0 & {\mbox{otherwise}}\:. \end{array} \right. \]
\begin{Prp} \label{prp1}
For all~$\tau, \tau' \in [\tau_0, \tau_1]$ with~$\tau' < \tau$,
the following inequality holds:
\[ \left( \int_{\tau'}^\tau A(\sigma)\, d\sigma \right)^2 \;\leq\;
\left( V(\tau) - V(\tau') \right) \left(F(\tau) - F(\tau') \right) . \]
\end{Prp}
{\Proof} The function~$V(\tau)$ is clearly monotone increasing.
Hence there is a unique Borel measure~$\nu$ such that
(see~\cite[Section~I.4]{RS})
\beq \label{Vdec}
V(\tau) \;=\; \int_{\tau_0}^\tau d\nu\:.
\eeq
The Lebesgue decomposition theorem~\cite[Section~I.4]{RS} allows us
to decompose this measure with respect to the Lebesgue measure,
\beq \label{RNT}
d\nu \;=\; g(\tau)\; d\tau + d\nu_{\mbox{\scriptsize{sing}}}
\quad {\mbox{with}} \quad 
g \in L^1([\tau_0, 1], d\tau) \:.
\eeq
If~$\tau$ is a regular value, we can compute~$g$ by differentiation,
\[ g(\tau) \;=\; \frac{d}{d\tau}\, V(\tau) \;=\;
\frac{d}{d\tau} \int^\tau d\sigma \int_{L(\sigma)} \frac{1}{|D\phi|}\: d\mu_{L(\sigma)} \;=\; \int_{L(\tau)} \frac{1}{|D\phi|}\: d\mu_{L(\tau)} \:. \]
Integrating~(\ref{RNT}) from~$\tau'$ to~$\tau$ and using~(\ref{Vdec}),
we find that
\beq \label{Ves}
V(\tau) - V(\tau') \;\geq\; \int_{\tau'}^\tau g(\sigma)\: d\sigma\:.
\eeq

We introduce the function~$\Theta_{\mbox{\scriptsize{reg}}}$ by
\[ \Theta_{\mbox{\scriptsize{reg}}}(\tau) \;=\;
\left\{ \begin{array}{cl}
1 & {\mbox{if~$\tau$ is a regular value}} \\
0 & {\mbox{otherwise}}\:. \end{array} \right. \]
Then the Schwarz inequality yields
\begin{eqnarray*}
\lefteqn{ \int_{\tau'}^\tau A(\sigma)\, d\sigma \;=\;
\int_{\tau'}^\tau \Theta_{\mbox{\scriptsize{reg}}}(\sigma) \:d\sigma 
\int_{L(\sigma)} d\mu_{L(\sigma)} } \\
&\leq&
\left( \int_{\tau'}^\tau \Theta_{\mbox{\scriptsize{reg}}}(\sigma) \:d\sigma
\int_{L(\sigma)} |D \phi|\:d\mu_{L(\sigma)} \right)^{\frac{1}{2}}
\left( \int_{\tau'}^\tau \Theta_{\mbox{\scriptsize{reg}}}(\sigma) \:d\sigma
\int_{L(\sigma)} \frac{1}{|D \phi|}\:
d\mu_{L(\sigma)} \right)^{\frac{1}{2}} \\
&=& \left( \int_{\tau'}^\tau f(\sigma)\, d\sigma \right)^{\frac{1}{2}}
\left( \int_{\tau'}^\tau g(\sigma)\, d\sigma \right)^{\frac{1}{2}} \:.
\end{eqnarray*}
Taking the square and using~(\ref{coarea}) as well as~(\ref{Ves})
gives the result.
\QED

\section{Applications to Curvature Estimates} \label{sec3}
For clarity, we begin with the simpler and physically interesting
case of dimension~$n=3$, the extension to higher dimension will be
given afterwards. Our starting point is the following integral estimate as
derived in~\cite{BF}.
\begin{Lemma} \label{lemma31}
Suppose~$\psi$ is a Witten spinor~(\ref{bvp}) in
a complete, asymptotically flat manifold~$(M^3, g)$ whose
scalar curvature is non-negative and integrable. Then
\[ \int_M |R|^2\; |\psi|^2\:d\mu_M \;\leq\;
m \:c_1\: \sup_M |R| \:+\:
\sqrt{m}\:c_2\: \|\nabla R\|_{L^2(M)} \;, \]
where the constants~$c_1$
and~$c_2$ are independent of the geometry of~$M$.
\end{Lemma}

\noindent
{\em{Proof of Theorem~\ref{thm31}. }}
If~$t_0<\tau_0$, we set~$t=t_0$. Then~$L(t)$ is empty and~$A(t)=0$.
If conversely~$t_0 \geq \tau_0$,
we apply Proposition~\ref{prp1} with~$\tau'=t_0$ and~$\tau=t_1$.
According to the mean value theorem, there is a subset of~$[t_0, t_1]$
of positive Lebesgue measure on which
\[ (t_1-t_0)\: A \;\leq\; \sqrt{(V(t_1) - V(t_0))\, (F(t_1)-F(t_0))}\:. \]
Out of this subset we choose a regular value~$t$. This gives~(\ref{Aineq}).
On the set~$M \setminus \Omega(t)$, the Witten spinor clearly
satisfies the bound~$|\psi| \geq t$. Using this in Lemma~\ref{lemma31}
gives the result.
\QED
Let us now prove the extension of Theorem~\ref{thm31} to
higher dimension, Theorem~\ref{thmg1}.
If~$n \geq 4$, the statement of Lemma~\ref{lemma31} no longer holds.
Instead, we must work with the spinor operator~$P_x$ as introduced in~\cite[Section~4]{FK}: We choose an orthonormal basis of constant spinors $(\psi^i_0)_{i=1,\ldots, N}$, $N=2^{[n/2]}$,
$\langle \psi_0^i, \psi_0^j \rangle = \delta^{ij}$,
and denote the corresponding solutions of the boundary problem~(\ref{bvp})
by $(\psi^i)_{i=1,\ldots,N}$. We define the {\em{spinor operator}} $P_x$ by
\beq \label{sodef}
P_x \;:\; S_xM \longrightarrow S_xM \;:\; \psi \longmapsto
    \sum_{i=1}^N \langle \psi^i_x,\: \psi \rangle\: \psi^i_x \;.
\eeq
Clearly, this operator is non-negative. We set
\beq \label{pdef}
p(x) \;=\; \inf \Big\{ \langle \chi, P_x \chi \rangle \:\Big|\:
\chi \in S_xM,\; |\chi|=1 \Big\} \;\geq\; 0 \:.
\eeq
Our starting point is the following curvature estimate, which is
an improvement of the estimates in~\cite{FK}.
\begin{Lemma} \label{lemmag1}
Suppose that~$(M^n, g)$, $n \geq 4$, is
a complete, asymptotically flat manifold whose scalar curvature
is non-negative and integrable.
Then the Riemann tensor~$R$ and the infimum of the
spinor operator~(\ref{sodef}, \ref{pdef}) satisfy the inequality
\[ \int_M |R|^2\; p(x)\:d\mu_M \;\leq\;
m \:c_1(n)\: \sup_M |R| \:+\:
\sqrt{m}\:c_2(n)\: \|\nabla R\|_{L^2(M)} \;, \]
where the constants~$c_1$ and~$c_2$ depend on the dimension,
but are independent of the geometry of~$M$.
\end{Lemma}
{\Proof} In~\cite[Corollary~3.2, Lemma~5.1]{FK} it was proved that,
choosing an orthonormal frame~$(s_\alpha)_{\alpha=1,\ldots, n}$,
\[ \int_M \sum_{\alpha, \beta=1}^n
\Tr \left( R^S(s_\alpha, s_\beta)^2\: P(x) \right) d\mu \;\leq\;
m \:c_1(n)\: \sup_M|R| \:+\:
\sqrt{m}\:c_2(n)\: \|\nabla R\|_{L^2(M)} \;, \]
where~$R^S$ is the curvature of the spin connection, which is related
to the Riemann tensor by
\[ R^S(X,Y)\: \psi=\frac{1}{4} \sum_{\alpha,\beta=1}^n R(X,Y,s_\alpha,s_\beta)
s_\alpha \cdot s_\beta \cdot \psi\:. \]
Introducing the abbreviation
\[ {\mathcal{R}}^2 \;=\; \sum_{\alpha, \beta=1}^n R^S(s_\alpha, s_\beta)^2 \:,\]
the operator~${\mathcal{R}}^2(x)$ acts on~$S_x(M)$ as a positive operator.
Its trace is a positive multiple of the norm squared
of the Riemann tensor,
\[ \Tr \left({\mathcal{R}}^2 \right) \;=\; c(n)\: |R|^2\:. \]
Hence
\[ \Tr \left({\mathcal{R}}^2\: P(x) \right) \;=\;
c(n)\: |R|^2\:p(x) \:+\: \Tr \left({\mathcal{R}}^2\: (P(x)-p(x)) \right)
\;\geq\; c(n)\: |R|^2\:p(x) \:, \]
where in the last step we used that the trace of the product of two
positive operators is positive.
\QED

\begin{Thm} \label{thmg1}
Let~$(M^n, g)$, $n \geq 4$,
be a complete, asymptotically flat manifold whose
scalar curvature is non-negative and integrable.
Suppose that for an interval~$[t_0, t_1]
\subset (0,1]$ there is a constant~$C$ such that every
Witten spinor~(\ref{bvp}) satisfies the volume bound
\[ V(t_1) - V(t_0) \;\leq\; C \:. \]
Then there is an open set~$\Omega \subset M$ with the following
properties. The $(n-1)$-dimensional Hausdorff measure~$\mu_{n-1}$ of the
boundary of~$\Omega$ is bounded by
\[ \mu_{n-1}(\partial \Omega) \;\leq\; \sqrt{m}\;c_0(n,t_0)\; \frac{\sqrt{C}}{t_1-t_0} \:. \]
On the set~$M \setminus \Omega$, the Riemann tensor satisfies the inequality
\[ \int_{M \setminus \Omega} |R|^2 \;\leq\;
m \:c_1(n,t_0)\: \sup_M |R| \:+\:
\sqrt{m}\:c_2(n,t_0)\: \|\nabla R\|_{L^2(M)} \:. \]
Here the constants~$c_0$, $c_1$ and~$c_2$ depend on the dimension and
on~$t_0$, but they are independent of the geometry of~$M$.
\end{Thm}
{\Proof} For given~$x \in M$ we introduce the mapping
\[ B \;:\; \C^N \rightarrow S_xM \;:\; \zeta \mapsto
\sum_{i=1}^N \zeta_i\: \psi_i(x)\:. \]
The a-priori bound~$|\psi(x)| \leq 1$ for all Witten spinors (see
the argument before~(\ref{22a})) yields that~$\|B\| \leq 1$.
Furthermore, the spinor operator can be written as~$P_x = B B^*$.
Since the operators~$B B^*$ and~$B^* B$ are both Hermitian and have
the same spectrum, we find
\[ p(x) \;=\; \inf \spec(B B^*) \;=\; \inf \spec(B^* B) \;=\;
\inf_{\zeta {\mbox{\scriptsize{ with }}} |\zeta|=1} |B \zeta|^2\:. \]

We choose a finite number of points~$\zeta^1, \ldots, \zeta^L$ on the
unit sphere~$S^N_1$ in~$\C^N$ such that the
balls $B_{t_0/2}(\zeta^1), \ldots, B_{t_0/2}(\zeta^L)$ cover~$S^N_1$
(with a constant~$L=L(n,t_0)$).
For every~$a \in \{1, \ldots, L\}$, we let~$\psi^a$ be the Witten
spinor~(\ref{bvp}) with boundary conditions~$\psi_0=\zeta^a$. Then
$\psi^a = B \zeta^a$.
Exactly as in the proof of Theorem~\ref{thm31}, for every~$\psi^a$ we
can choose a regular value~$t \in [t_0, t_1]$ such that
\[ A(t) \;\leq\; \frac{\sqrt{m C}}{t_1-t_0} \:. \]
We also denote~$\Omega(t)$ by~$\Omega^a$ and set~$\Omega= \cup_{a=1}^L
\Omega^a$. Then~$\partial \Omega$ is a subset of~$\cup_{a=1}^L \partial
\Omega^a$, and thus its Hausdorff measure is bounded by
\[ \mu_{n-1}(\partial \Omega) \;\leq\; L\, \frac{\sqrt{m C}}{t_1-t_0} \:. \]

In view of the estimates of Lemma~\ref{lemmag1}, it remains to show that
\[ p(x) \;\geq\; \frac{t_0^2}{4} \spc \forall x \in M \setminus \Omega\:. \]
For any~$\zeta \in S^N_1$ we can choose an index~$a \in \{1, \ldots, L\}$
such that~$|\zeta - \zeta^a| < t_0/2$. Thus
\begin{eqnarray*}
|B \zeta| &\geq& |B \zeta^a| - | B (\zeta- \zeta^a)|
\;=\; |\psi^a| - | B (\zeta- \zeta^a)| \\
&\geq& |\psi^a| - \| B\|\,|\zeta- \zeta^a|
\;\geq\; t_0 - \| B\|\,|\zeta- \zeta^a| \;\geq\; \frac{t_0}{2}\:,
\end{eqnarray*}
where in the last step we used that~$\|B\| \leq 1$.
\QED

\section{Estimates in a Conformal Compactification} \label{sec4}
This section is devoted to the proof of Theorem~\ref{thm32} and its
generalization to higher dimension, Theorem~\ref{thmg2}.
By rescaling, we can arrange that the total mass~$m$ equals two.
We assume that~$(M^n, g)$, $n \geq 3$, is asymptotically
Schwarzschild, i.e.\ there is a diffeomorphism
\[ \varphi : M \setminus K \to \R^n \setminus \overline{B_{\rho}(0)} \]
such that
\begin{equation} \label{Schmet}
(\varphi_* g)_{ij} \;=\; \left(1+\frac{1}{|x|^{n-2}}\right)^{\frac{4}{n-2}}\:
\delta_{ij}\:.
\end{equation}
On~$M \setminus K$ we introduce the function~$r(x) = |\phi(x)|$.
For the point compactification, we choose parameters~$\sigma$ and~$R$ in the
range
\beq \label{scaling}
\rho \;\leq\; \sigma \;\leq \; R \;<\; c(n)\, \rho
\eeq
and consider a function~$\lambda$ with the following properties
(for the construction of~$\lambda$ see~\cite{FKr}):
\begin{description}
\item[(i)] $\displaystyle \lambda_{|K} \;\equiv\; 1$
\item[(ii)] $\displaystyle
\lambda(x) \;=\; \left(\frac{2 \sigma^2}{\sigma^2+r(x)^2}\right)\cdot
 \left(1+\frac{1}{r(x)^{n-2}}\right)^{-\frac{2}{n-2}}
\qquad {\mbox{on~$\phi^{-1}(\R^n \setminus B_{R}(0))$}}.$
\item[(iii)] The scalar curvature corresponding to the conformally
changed metric
\beq \label{cc}
\tilde{g} \;=\; \lambda^2\: g
\eeq
is non-negative.
\end{description}
After the conformal change, the region~$r>R$ is isometric to
a neighborhood of the north pole~$\mathfrak{n}$
of the sphere~$S^n_\sigma$ with the north pole removed.
Adding the north pole, we obtain the conformal
compactification~$(\bar{M}, \tilde{g})$.
For any~$r \geq R$ we introduce the spherical cap
\[ C_r \;=\; \ov{(\varphi^{-1}(\R^n \setminus B_r(0)),
\tilde{g})} \subset \bar M \:. \]
Finally, we denote the geodesic radius of the spherical
cap~$C_R \subset S^n_\sigma$ by~$\delta$.

Our main task is to bound the function~$\phi$ in the asymptotic end from below,
see Proposition~\ref{prplower}.
As in~\cite{FKr} we work on the sphere~$S^n_\sigma$ with
Sobolev norms which are scaling invariant in~$\sigma$, namely
\[ \|f\|^2_{H^{k,2}(S^n_\sigma)} \;:=\; \sum_{\kappa {\mbox{\scriptsize{ with }}} |\kappa|\leq k} \sigma^{2 |\kappa| - n}
\int_{S^n_\sigma} \|\nabla^\kappa f(x)\|^2 \:dx \:. \]
We let~$\eta$ be a smooth function on~$\bar{M}$ with
\[ {\mbox{supp}}\, \eta \subset C_R \:,\spc
\eta |_{C_{2R}} \;\equiv\; 1\:. \]
The next lemma is an elliptic estimate for the Dirac operator
on~$S^n_\sigma$, for the proof see~\cite[Lemma~5.1]{FKr}.
\begin{Lemma} \label{lemmapower}
For any smooth section~$\chi$ in~$S(S^n_\sigma)$,
\[ \|\eta^{k+1} \chi \|^2_{H^{k,2}(S^n_\sigma)} \;\leq\;
c(n) \sum_{l=0}^k \sigma^{2 l - n} \;\|\eta^{l+1}\:
\tilde{\D}^l \chi \|_{L^2(S^n_\sigma)}^2 \:. \]
\end{Lemma}
For a given Witten spinor~$\psi$ we introduce the spinors
\begin{eqnarray}
\psi_{\mbox{\scriptsize{asy}}} &=& \eta \left( 1 + \frac{1}{r(x)^{n-2}}
\right)^{-\frac{n-1}{n-2}} \psi_0 \spc\, {\mbox{on $M \setminus K$}}
\label{asydef} \\
\delta \psi &=& \psi - \eta\, \psi_{\mbox{\scriptsize{asy}}} 
\spc\spc\spc\quad {\mbox{on $M$}}\:, \label{dpdef}
\end{eqnarray}
where~$\psi_0$ is the boundary value of~$\psi$ at infinity.
The spinor~$\delta \psi$ is referred to as the {\em{Witten deviation}}.
We shall also consider the above spinors on the conformal
compactification~$(\bar{M}, \tilde{g})$. We then denote them with
an additional tilde and rescale them as usual by
\beq \label{rescale}
\tilde{\psi} \;=\; \lambda^{\frac{1-n}{2}}\: \psi\:.
\eeq
Since the metric on~$M \setminus K$ is Schwarzschild, the
conformal invariance of the Dirac equation shows
that~$\D \psi_{\mbox{\scriptsize{asy}}} = 0$,
and therefore also~$\tilde{\D} \tilde{\psi}_{\mbox{\scriptsize{asy}}} = 0$.
As a consequence,
\[ \tilde{\D} (\widetilde{\delta \psi}) \;=\;
-(\tilde{\D} \eta)\, \tilde{\psi}_{\mbox{\scriptsize{asy}}} \;=:\;
h\:. \]
It is important that the function~$h$ appearing here can be given
explicitly and is supported inside the spherical
cap~$C_R$.
\begin{Lemma} \label{lemmacap}
There is a constant~$c$ depending only on~$n$ and
the ratio~$\delta/\sigma$ such that
\[ \sup_{C_{2R}} \left( \lambda^{\frac{1-n}{2}}\,
|\delta \psi| \right) \;\leq\; \frac{c}{
\sigma\: \inf \spec |\tilde{\D}|} \:. \]
\end{Lemma}
{\Proof} In view of~(\ref{rescale}), we must estimate the
sup-norm of~$\widetilde{\delta \psi}$ on~$C_{2R}$. Obviously,
\[ \sup_{C_{2R}} |\widetilde{\delta \psi}| \;\leq\;
\sup_{C_R} |\eta^{k+1} \widetilde{\delta \psi}| \:. \]
Extending the last function by zero to~$S^n_\sigma$, we can apply
the Sobolev imbedding theorem on the sphere~$S^n_\sigma$. Thus
for sufficiently large~$k$,
\[ \sup_{C_R} |\eta^{k+1} \widetilde{\delta \psi}| \;=\;
\sup_{S^n_\sigma} |\eta^{k+1} \widetilde{\delta \psi}| \;\leq\;
c \,\| \eta^{k+1} \widetilde{\delta \psi} \|_{H^{k,2}(S^2_\sigma)}\:, \]
where~$c=c(n)$ is the Sobolev constant on the unit sphere.
Applying Lemma~\ref{lemmapower}, we obtain
\begin{eqnarray*}
\sup_{C_{2R}} |\widetilde{\delta \psi}|^2 &\leq& c^2
\sum_{l=0}^k \sigma^{2 l - n} \;\|\eta^{l+1}\:
\tilde{\D}^l (\widetilde{\delta \psi}) \|_{L^2(S^n_\sigma)}^2 \\
&=& c^2 \sigma^{- n} \;\|\eta\: (\widetilde{\delta \psi}) \|_{L^2(S^n_\sigma)}^2
\:+\: c^2 \sum_{l=1}^k \sigma^{2 l - n} \;\|\eta^{l+1}\:
\tilde{\D}^{l-1} h \|_{L^2(S^n_\sigma)}^2\:.
\end{eqnarray*}
The obtained terms can be estimated as follows,
\begin{eqnarray*}
\|\eta\: (\widetilde{\delta \psi}) \|_{L^2(S^n_\sigma)}^2 &\leq&
\| \widetilde{\delta \psi} \|_{L^2(M)}^2 \;\leq\;
\frac{1}{\inf \spec(\tilde{\D}^2)}\: \| h \|_{L^2(S^n_\sigma)}^2
\;\leq\; \frac{c\, \sigma^{n-2}}{\inf \spec(\tilde{\D}^2)} \\
\|\eta^{l+1}\: \tilde{\D}^{l-1} h \|_{L^2(S^n_\sigma)}^2 &\leq&
\|\tilde{\D}^{l-1} h \|_{L^2(S^n_\sigma)}^2 \;\leq\;
c\, \sigma^{n-2l}\:.
\end{eqnarray*}
Putting these estimates together and taking the square root
gives the estimate
\beq \label{prelim}
\sup_{C_{2R}} \left( \lambda^{\frac{1-n}{2}}\,
|\delta \psi| \right) \;\leq\; c \left[ 1 + \frac{1}{
\sigma\: \inf \spec |\tilde{\D}|} \right] \:.
\eeq
Finally, we can use the lower spectral bound
(see~\cite[Proof of Theorem~7.5]{FKr})
\beq \label{lsb}
\inf \spec (\tilde{D}^2) \;\leq\; \frac{c}{\sigma^2}
\eeq
to drop the first term in the square brackets in~(\ref{prelim}).
\QED

\begin{Prp} \label{prplower} There is a constant~$c$ depending only
on~$n$ and the quotient~$\delta/\sigma$ such that for all~$x \in M \setminus K$
with
\[ r(x) \;\geq\; r_1 \::=\: c \,\sigma \left(\sigma \inf \spec |\tilde{D}| \right)^{-\frac{1}{n-1}} \:, \]
the norm of the Witten spinor is bounded from below by
\[ \phi(x) \;\geq\; \frac{1}{2} \:. \]
\end{Prp}
{\Proof} We estimate~$\phi$ from below by
\[ \phi(x) \;=\; |\psi(x)| \;\geq\;
|\psi_{\mbox{\scriptsize{asy}}}| - |\delta \psi|\:. \]
As is obvious from the definition of~$\psi_{\mbox{\scriptsize{asy}}}$, (\ref{asydef}), by choosing~$c$ sufficiently large we can arrange
that~$|\psi_{\mbox{\scriptsize{asy}}}| > 3/4$.
Hence we need to arrange that~$|\delta \psi|<1/4$.
According to Lemma~\ref{lemmacap}, this can be achieved by choosing
\[ \lambda^{\frac{1-n}{2}} \;\geq\; \frac{4 c}
{\sigma\: \inf \spec |\tilde{\D}|}\: . \]
Using the explicit form of~$\lambda$ in~{\bf{(ii)}} gives the result.
\QED
Next we derive the desired volume bound.
\begin{Lemma} \label{lemmavb}
There is a constant~$c$ which depends only on~$n$ but is independent
of the geometry of~$M$ such that
\[ V \Big( \frac{1}{2} \Big) - V \Big( \frac{1}{4} \Big) \;\leq\;
c\: \frac{(\rho+1)^n}{\sigma^2\, \inf \spec ({\mathcal{D}}^2)} \:. \]
\end{Lemma}
{\Proof} Since~$|\psi| \geq 1/4$ on~$\Omega(1/2) \setminus \Omega(1/4)$,
we clearly have
\[ V \Big( \frac{1}{2} \Big) - V \Big( \frac{1}{4} \Big) \;\leq\; 16 \int_{\Omega(1/2)} |\psi|^2\: d\mu_M \:. \]
Furthermore, applying Proposition~\ref{prplower},
\[ \int_{\Omega(1/2)} |\psi|^2\: d\mu_M \;\leq\;
\int_{M \setminus C_{r_1}} |\psi|^2\: d\mu_M \;\leq\;
\int_K |\psi|^2\: d\mu_M \:+\: \mu \Big(
\{x \in M \setminus K \;|\; r(x) \leq r_1 \} \Big) , \]
where in the last step we used that~$|\psi| \leq 1$, (\ref{22a}).
To the integral over~$K$ we apply the weighted $L^2$-estimates
in~\cite[Corollary~7.6]{FKr}, whereas the additional measure can be
estimated by the volume of a Euclidean ball in the asymptotic end,
\begin{eqnarray*}
&& \int_K |\psi|^2\: d\mu_M \;\leq\; \int_K |\psi|^2\: d\mu_M
\:+\: \int_{M \setminus K} |\delta \psi|^2\: d\mu_M
\;\leq\; c\: \frac{(\rho+1)^n}{\sigma^2\, \inf \spec (\tilde{D}^2)} \\
&& \mu \Big( \{x \in M \setminus K \;|\; r(x) \leq r_1 \} \Big)
\;\leq\; c\, r_1^n \;\leq\; c\, (\rho+1)^n
\left(\sigma \inf \spec |\tilde{D}| \right)^{-\frac{n}{n-1}}\:,
\end{eqnarray*}
where in the last step we used that, according to~(\ref{scaling}),
$\sigma$ and~$\rho$ have the same scaling. Combining these inequalities
with~(\ref{lsb}), we obtain the result.
\QED

\vspace*{.3em} \noindent
{\em{Proof of Theorem~\ref{thm32}. }}
We apply Theorem~\ref{thm31} with~$t_0=1/4$ and~$t_1=1/2$, using the
estimate~$F(t_1) - F(t_0) \leq m$. We then put in the
estimate of Lemma~\ref{lemmavb}.
\QED
Applying Lemma~\ref{lemmavb} in the same way to Theorem~\ref{thmg1}
gives the following result.

\begin{Thm} \label{thmg2}
Let~$(M^n, g)$, $n \geq 4$, be a complete manifold of non-negative
scalar curvature such that~$M \setminus K$ is isometric to
the Schwarzschild geometry.
Then there is an open set~$\Omega \subset M$ with the following
properties. The $(n-1)$-dimensional Hausdorff measure~$\mu_{n-1}$ of the
boundary of~$\Omega$ is bounded by
\[ \mu_{n-1}(\partial \Omega) \;\leq\; \;c_0(n)\: \sqrt{m}\;
\frac{\left( \rho+m^{\frac{1}{n-2}} \right)^{\frac{n}{2}}}
{\sigma\, \inf \spec |\tilde{\D}|} \:. \]
On the set~$M \setminus \Omega$, the Riemann tensor satisfies the inequality
\[ \int_{M \setminus \Omega} |R|^2 \;\leq\;
m \:c_1(n)\: \sup_M |R| \:+\:
\sqrt{m}\:c_2(n)\: \|\nabla R\|_{L^2(M)} \:. \]
Here the constants~$c_0$, $c_1$ and~$c_2$ depend on the dimension,
but are independent of the geometry of~$M$.
\end{Thm} 

%#################################################################

\noindent
NWF I -- Mathematik,
Universit{\"a}t Regensburg, 93040 Regensburg, Germany, \\
{\tt{Felix.Finster@mathematik.uni-regensburg.de}}


\begin{thebibliography}{99}
\bibitem{ADM} R.\ Arnowitt, S.\ Deser, C.\ Misner, ``Energy and the
Criteria for Radiation in General Relativity,'' {\em{Phys.\ Rev.}}\ 118,
1100 (1960)
\bibitem{B} R.\ Bartnik, ``The Mass of an
Asymptotically Flat Manifold,'' {\em{Commun.\ Pure Appl.\ Math.}}~XXXIX
(1986) 661--693
\bibitem{BF} H.\ Bray, F.\ Finster, ``Curvature estimates and the
positive mass theorem,'' math.DG/9906047,
{\em{Comm.\ Anal.\ Geom.}}\ {\bf{10}} (2002) 291--306
\bibitem{FK} F.\ Finster, I.\ Kath, ``Curvature estimates in asymptotically
flat manifolds of positive scalar curvature,'' math.DG/0101084,
{\em{Comm.\ Anal.\ Geom.}}\ {\bf{10}} (2002) 1017--1031
\bibitem{FKr} F.\ Finster, M.\ Kraus, ``A weighted $L^2$-estimate of the Witten spinor in asymptotically Schwarzschild manifolds,'' math.DG/0501195, 
{\em{Canadian J.\ Math.}} {\bf{59}} (2007) 943--965
\bibitem{HI} G.\ Huisken, T.\ Ilmanen, ``The inverse mean curvature flow and the Riemannian Penrose inequality,'' {\em{J.\ Differential Geom.}}\ {\bf{59}} (2001) 353--437
\bibitem{LM} H.-B. Lawson, M.-L. Michelsohn, ``Spin Geometry,''
{\em{Princeton University Press}}, Princeton (1989)
\bibitem{L} J.\ Lohkamp, ``The higher dimensional positive mass theorem I,''
math.DG/0608795 (2006)
\bibitem{PT} T.\ Parker, C.\ H.\ Taubes, ``On Witten's proof of the
positive energy theorem,'' {\em{Comm.\ Math.\ Phys.}}\ {\bf{84}} (1982) 223--238
\bibitem{RS} M.\ Reed, B.\ Simon, ``Methods of Modern Mathematical
Physics, I: Functional Ana\-lysis,'' {\em{Academic Press}} (1980)
\bibitem{SY} R.\ Schoen, S.-T.\ Yau, ``On the proof of the positive
mass conjecture in general relativity,'' {\em{Comm.\ Math.\ Phys.}}\
{\bf{65}} (1979) 45--76
\bibitem{W} E.\ Witten, ``A new proof of the positive energy
theorem,'' {\em{Comm.\ Math.\ Phys.}}\ {\bf{80}} (1981) 381--402
\end{thebibliography}
\end{document}